\documentclass[12pt,reqno]{amsart}

\usepackage{amssymb}

\usepackage{eucal}

\input{diagrams}

\def\sD{\mbox{\sf D}}
\def\sE{\mbox{\sf E}}

\def\sG{\mbox{\sf G}}

\def\sK{\mbox{\sf K}}

\def\Ab{\mbox{\sf Ab}}

\def\ast{{\textstyle *}}

\def\D{\sD}

\def\H{\operatorname{H}}

\def\Hom{\operatorname{Hom}}
\def\id{\operatorname{id}}

\def\Ker{\operatorname{Ker}}

\def\LTensor{\stackrel{\operatorname{L}}{\otimes}}
\def\Mod{\mbox{\sf Mod}}

\def\opp{\operatorname{op}}

\def\Proj{\mbox{\sf Pro}}

\def\RHom{\operatorname{RHom}}

\numberwithin{equation}{part}

\newtheorem{Lemma}{Lemma}[section]
\newtheorem{Theorem}[Lemma]{Theorem}
\newtheorem{Proposition}[Lemma]{Proposition}
\newtheorem{Corollary}[Lemma]{Corollary}

\theoremstyle{definition}
\newtheorem{Definition}[Lemma]{Definition}
\newtheorem{Setup}[Lemma]{Setup}

\newtheorem{Remark}[Lemma]{Remark}

\def\generalR{R}

\def\R{A}
\def\matlisR{B}

\def\genericcategory{\sK}
\def\genericsubcategory{\sE}

\def\KProjR{\sK(\Proj\,\R)}
\def\EnochsR{\sE(\R)}

\def\inc{e_{\ast}}
\def\adj{e^!}

\def\Rm{$\R$-module}
\def\Rlm{$\R$-left-module}

\def\projmod{Q}         
\def\injmod{J}          
\def\projcomplx{P}      
\def\projcomplxtil{\widetilde{P}}
\def\injcomplx{I}       
\def\Enochscomplx{E}    
\def\Enochscomplxtil{\widetilde{E}}      
\def\Gorprojmod{G}      
\def\Gorprojmodtil{\widetilde{G}}        
\def\flatmod{F}         

\def\gtil{\widetilde{g}}
\def\etil{\widetilde{e}}
\def\gammatil{\widetilde{\gamma}}
\def\epsilontil{\widetilde{\epsilon}}

\begin{document}

\title[Gorenstein projective modules]
{The Gorenstein projective modules are precovering}

\author{Peter J\o rgensen}
\address{Department of Pure Mathematics, University of Leeds,
Leeds LS2 9JT, United Kingdom}
\email{popjoerg@maths.leeds.ac.uk, www.maths.leeds.ac.uk/\~{ }popjoerg}




\begin{abstract} 
The Gorenstein projective modules are proved to form a precovering
class in the module category of a ring which has a dualizing complex.
\end{abstract}

\maketitle

\setcounter{section}{-1}
\section{Introduction}
\label{sec:introduction}

This paper proves over a wide class of rings that the Gorenstein
projective modules form a precovering class in the module category.
Let me explain this statement.  There are two terms of mystery,
``Gorenstein projective modules'' and ``precovering class''; I will
explain the latter first.

{\em Precovering classes} are also known by the name of contravariantly
finite classes.  In a module category, a class $\sG$ of modules is
precovering if is satisfies the following: For each module $M$, there
exists a homomorphism $\Gorprojmod \longrightarrow M$ with
$\Gorprojmod$ in $\sG$, such that if $\Gorprojmodtil \longrightarrow M$
is any homomorphism with $\Gorprojmodtil$ in $\sG$ then the dotted
arrow exists to make the following diagram commutative,
\[
  \begin{diagram}
                   &           & \Gorprojmod \\
                   & \ruDotsto & \dTo \\
    \Gorprojmodtil & \rTo      & M \lefteqn{.}
  \end{diagram}
\smallskip
\]
The homomorphism $\Gorprojmod \longrightarrow M$ is called a
$\sG$-precover of $M$.  By taking the kernel of $\Gorprojmod
\longrightarrow M$, taking a $\sG$-precover, and repeating, I can
construct a $\sG$-resolution of $M$,
\[
  \cdots 
  \longrightarrow \Gorprojmod_2 
  \longrightarrow \Gorprojmod_1 
  \longrightarrow \Gorprojmod_0
  \longrightarrow M 
  \rightarrow 0,
\]
which becomes exact when I apply the functor $\Hom(\Gorprojmod,-)$
with $\Gorprojmod$ in $\sG$.  The name ``precovering class'' is due to
\cite{Enochs}.  If $\sG$ is precovering, then I can do homological
algebra using $\sG$-resolutions instead of projective resolutions.
There is a large literature on this so-called ``relative homological
algebra''; an example is the recent \cite{HHGorensteinDerived}.

{\em Gorenstein projective modules} are modules which have the form
$\Gorprojmod = \Ker(\Enochscomplx^1 \longrightarrow \Enochscomplx^2)$,
where $\Enochscomplx$ is a complex of projective modules which is
exact and satisfies that $\Hom(\Enochscomplx,\projmod)$ is exact for
each projective module $\projmod$.  The first complete statement of
this definition seems to be in \cite{EJX}, but the idea goes back to
\cite{AuslanderBridger}.  The point is that one can do homological
algebra with Gorenstein projective modules instead of projective
modules, and that such ``Gorenstein homological algebra'' does for
Gorenstein rings what ordinary homological algebra does for rings of
finite global dimension.  The archetypal result is that a noetherian
local commutative ring is Gorenstein if and only if each finitely
generated module has a bounded Gorenstein projective resolution.  An
extensive theory of Gorenstein homological algebra has been developed;
part of it is already in \cite{AuslanderBridger}, a comprehensive
source which was up to date when it was written is \cite{Winther}, and
for some recent work see \cite{HHGorensteinDerived} and
\cite{HHGorensteinHomDim}. 

A weakness of the existing literature on Gorenstein homological
algebra is that it has only been known in special cases that the
Gorenstein projective modules form a precovering class.  The state of
the art appears to be \cite[prop.\ 2.18]{HHGorensteinHomDim} which
only proves the precovering property over Gorenstein rings.

The traditional remedy for this weakness has been to relax the
conditions imposed on the Gorenstein projective resolutions used in
the theory, which are then required just to be exact rather than
obtained from successive precovers.  This is formalized in the theory
of resolving classes, but suffers from the serious shortcoming that it
does not permit the definition of relative derived functors.

However, this paper removes the weakness by proving that the
Gorenstein projective modules do form a precovering class over a ring
which has a dualizing complex.  For the sake of simplicity, the main
part of the paper, sections \ref{sec:lemma} to \ref{sec:modules},
proves this result over a noetherian commutative ring.  However, as I
will show in section \ref{sec:noncomm}, the proofs really apply to
much more general (non-commutative) algebras with dualizing complexes.

The idea of the proof is taken from \cite{PJSpectra}, and is, to my
knowledge, different from that used in other papers on precovering
classes.  Rather than attack the problem directly, I pass to
$\KProjR$, the homotopy category of complexes of projective modules
over the ring $\R$.

Inside it sits the subcategory $\EnochsR$ of complexes $\Enochscomplx$
which are exact and have $\Hom(\Enochscomplx,\projmod)$ exact for each
projective module $\projmod$.  Crucially, this subcategory can be
characterized as the kernel of a homological functor which respects
small coproducts (see the proof of proposition \ref{pro:adjoint}), and
this enables me to use Bousfield localization to see that the
inclusion functor $\EnochsR \longrightarrow \KProjR$ has a right
adjoint $\KProjR \longrightarrow \EnochsR$ (proposition
\ref{pro:adjoint}).

This implies that $\EnochsR$ is a precovering class in $\KProjR$
(proposition \ref{pro:Enochs_is_precovering}).  However, the
Gorenstein projective modules are the modules of the form
$\Ker(\Enochscomplx^1 \longrightarrow \Enochscomplx^2)$ for
$\Enochscomplx$ in $\EnochsR$, and it turns out that the result on
$\EnochsR$ descends to give that the Gorenstein projective modules
form a precovering class in the module category (lemma
\ref{lem:Gorenstein_projectives_are_precovering} and theorem
\ref{thm:Gorenstein_projectives_are_precovering}).

Let me mention some related work: First, several previous papers have
investigated whether the Gorenstein projective modules form a
precovering class.  As mentioned, I believe the state of the art to be
\cite[prop.\ 2.18]{HHGorensteinHomDim}.  Secondly, in addition to
defining Gorenstein projective modules, the paper \cite{EJX} also
defined Gorenstein flat and Gorenstein injective modules, and
\cite{EnochsJendaLopezRamos} and \cite{EnochsLopezRamos} proved
for large classes of rings that the Gorenstein flat modules form a
precovering class and that the Gorenstein injective modules form a
preenveloping class (the dual notion to precovering).  Hence the
present paper is a natural complement to \cite{EnochsJendaLopezRamos}
and \cite{EnochsLopezRamos}.

\medskip
\noindent
{\em Notation. } 
Let me close the introduction by setting up a minimum of notation.  In
sections \ref{sec:lemma}, \ref{sec:complexes}, and
\ref{sec:modules} (but not in section \ref{sec:noncomm}), the
following two setups are in force.

\begin{Setup}
\label{set:blanket1}
Let $\R$ be a noetherian commutative ring with a dualizing complex
$D$.  That is,
\begin{enumerate}

  \item  The cohomology of $D$ is bounded and finitely generated
         over $\R$.

  \item  The injective dimension $\id_{\R} D$ is finite.

  \item  The canonical morphism $\R \longrightarrow
         \RHom_{\R}(D,D)$ in the derived category $\D(\R)$ is an
         isomorphism. 

\end{enumerate}
\end{Setup}

\begin{Setup}
\label{set:I1}
Let $D \stackrel{\simeq}{\longrightarrow} \injcomplx$ be an injective
resolution so that $\injcomplx$ is a bounded complex.
\end{Setup}

See \cite[chp.\ V]{HartsResDual} for background on dualizing
complexes. 

\begin{Definition}
\label{def:Enochs}
By $\EnochsR$ is denoted the class of complexes $\Enochscomplx$ of \Rm
s so that $\Enochscomplx$ consists of projective modules, is exact,
and has $\Hom_{\R}(\Enochscomplx,\projmod)$ exact for each projective
\Rm\ $\projmod$.
\end{Definition}

I will view $\EnochsR$ as a full subcategory of $\KProjR$, the
homotopy category of complexes of projective \Rm s.

\begin{Definition}
\label{def:Gorenstein_projectives}
An \Rm\ is called Gorenstein projective if it has the form
$\Ker(\Enochscomplx^1 \longrightarrow \Enochscomplx^2)$ for some
$\Enochscomplx$ in $\EnochsR$.
\end{Definition}

Observe that each projective \Rm\ $\projmod$ is Gorenstein projective,
since $\projmod$ is equal to $\Ker(\projmod \longrightarrow 0)$, and
since $\projmod \longrightarrow 0$ is part of the complex $\cdots
\longrightarrow 0 \longrightarrow \projmod
\stackrel{\id}{\longrightarrow} \projmod \longrightarrow 0
\longrightarrow \cdots$ which is null homotopic and hence in
$\EnochsR$.

\begin{Remark}
\label{rmk:flats_have_finite_pd}
Since $\R$ has a dualizing complex, it has finite Krull dimension
by \cite[cor.\ V.5.2]{HartsResDual}, so by \cite[Seconde partie, cor.\
(3.2.7)]{RaynaudGruson}, each flat \Rm\ has finite projective
dimension.
\end{Remark}

\section{A lemma}
\label{sec:lemma}

The following lemma uses $\injcomplx$, the bounded injective
resolution of $D$ from setup \ref{set:I1}.

\begin{Lemma}
\label{lem:technical}
Let $\projcomplx$ be a complex of projective \Rm s.  Then
\begin{eqnarray*}
  & \mbox{$\Hom_{\R}(\projcomplx,\projmod)$ is exact for each
          projective \Rm\ $\projmod$} & \\
  & \mbox{$\Leftrightarrow \injcomplx \otimes_{\R} \projcomplx$
          is exact.} &
\end{eqnarray*}
\end{Lemma}

\begin{proof}
\noindent
$\Rightarrow \;$  Suppose that $\Hom(\projcomplx,\projmod)$ is exact
for each projective module $\projmod$.  To see that $\injcomplx
\otimes \projcomplx$ is an exact complex, it is enough to see that
\[
  \Hom(\injcomplx \otimes \projcomplx,\injmod) 
  \cong \Hom(\projcomplx,\Hom(\injcomplx,\injmod))
\]
is exact for each injective module $\injmod$.  

But $\Hom(\injcomplx,\injmod)$ is a bounded complex of flat modules,
so is finitely built from flat modules in the homotopy category of
complexes of \Rm s, $\sK(\Mod\,\R)$, so it is enough to see that
$\Hom(\projcomplx,\flatmod)$ is exact for each flat module $\flatmod$.

Since $\flatmod$ has finite projective dimension by remark
\ref{rmk:flats_have_finite_pd}, there is a projective resolution
$\projcomplxtil \stackrel{\simeq}{\longrightarrow} \flatmod$ with
$\projcomplxtil$ bounded.  Since $\projcomplx$ consists of projective
modules and both $\projcomplxtil$ and $\flatmod$ are bounded, this
induces a quasi-isomorphism
\[
  \Hom(\projcomplx,\projcomplxtil) \simeq \Hom(\projcomplx,\flatmod).
\]
So it is enough to see that $\Hom(\projcomplx,\projcomplxtil)$ is exact.  

But $\projcomplxtil$ is a bounded complex of projective modules, so is
finitely built from projective modules, so it is enough to see that
$\Hom(\projcomplx,\projmod)$ is exact for each projective module
$\projmod$.  And this holds by assumption.

\smallskip

\noindent
$\Leftarrow \;$ Suppose that $\injcomplx \otimes \projcomplx$ is an
exact complex.  I must show that $\Hom(\projcomplx,\projmod)$ is exact
for each projective module $\projmod$.

First observe that by \cite[thm.\ (3.2)]{AvrFoxPLMS}, there is an
isomorphism
\[
  \projmod 
  \stackrel{\sim}{\longrightarrow} \RHom(D,D \LTensor \projmod).
\]
Of course, I can replace $D$ by $\injcomplx$ to get
\[
  \projmod 
  \stackrel{\sim}{\longrightarrow} 
  \RHom(\injcomplx,\injcomplx \LTensor \projmod).
\]
Here $\injcomplx \LTensor \projmod \cong \injcomplx \otimes \projmod$
because $\projmod$ is projective.  Moreover, $\injcomplx \otimes
\projmod$ is a bounded complex of injective modules so
$\RHom(\injcomplx,\injcomplx \LTensor \projmod) \cong
\RHom(\injcomplx,\injcomplx \otimes \projmod) \cong
\Hom(\injcomplx,\injcomplx \otimes \projmod)$.  So the above
isomorphism in the derived category is represented by the chain map
\[
  \projmod 
  \longrightarrow
  \Hom(\injcomplx,\injcomplx \otimes \projmod)
\]
which must accordingly be a quasi-isomorphism.  

Completing to a distinguished triangle in $\sK(\Mod\,\R)$ gives
\[
  \projmod 
  \longrightarrow \Hom(\injcomplx,\injcomplx \otimes \projmod) 
  \longrightarrow C 
  \longrightarrow
\]
where $C$ is exact.  Here $\injcomplx$ and $\injcomplx \otimes
\projmod$ are bounded, so $\Hom(\injcomplx,\injcomplx \otimes
\projmod)$ is bounded.  As the same is true for $\projmod$, the
complex $C$ is also bounded.

Now, the distinguished triangle gives another distinguished triangle
\[
  \Hom(\projcomplx,\projmod) 
  \longrightarrow
  \Hom(\projcomplx,\Hom(\injcomplx,\injcomplx \otimes \projmod)) 
  \longrightarrow
  \Hom(\projcomplx,C) 
  \longrightarrow.
\]
Here $\Hom(\projcomplx,C)$ is exact because $\projcomplx$ is a complex
of projective modules while $C$ is a bounded exact complex.  So to see
that $\Hom(\projcomplx,\projmod)$ is exact as desired, it is enough to
see that $\Hom(\projcomplx,\Hom(\injcomplx,\injcomplx \otimes
\projmod))$ is exact.

However,
\[
  \Hom(\projcomplx,\Hom(\injcomplx,\injcomplx \otimes \projmod)) 
  \cong 
  \Hom(\injcomplx \otimes \projcomplx,\injcomplx \otimes \projmod).
\]
And this is exact because $\injcomplx \otimes \projcomplx$ is exact by
assumption while $\injcomplx \otimes \projmod$ is a bounded complex of
injective modules.
\end{proof}

\section{Complexes}
\label{sec:complexes}

\begin{Lemma}
\label{lem:cg}
The triangulated category $\KProjR$ is compactly generated.
\end{Lemma}

\begin{proof}
The ring $\R$ is noetherian and hence coherent, and by remark
\ref{rmk:flats_have_finite_pd} each flat \Rm\ has finite projective
dimension.  So the proposition holds by \cite[thm.\ 2.4]{PJhomproj}. 
\end{proof}

Combining lemmas \ref{lem:technical} and \ref{lem:cg} gives the
following key result.

\begin{Proposition}
\label{pro:adjoint}
The inclusion functor $\inc : \EnochsR \longrightarrow \KProjR$ has a
right-adjoint $\adj : \KProjR \longrightarrow \EnochsR$.
\end{Proposition}

\begin{proof}
Consider the functor
\[
  k(-) = \H^0((\R \oplus \injcomplx) \otimes_{\R} -) : 
  \KProjR \longrightarrow \Ab
\]
from the homotopy category of complexes of projective \Rm s to the
category of abelian groups.  This is clearly a homological functor
respecting set indexed coproducts.  Moreover,
\[
  k(\Sigma^i \projcomplx) 
  \cong 
  \H^i(\projcomplx) \oplus \H^i(\injcomplx \otimes \projcomplx),
\]
where $\Sigma^i$ denotes $i$'th suspension, so for $\projcomplx$ to
satisfy $k(\Sigma^i \projcomplx) = 0$ for each $i$ means
\[
  \H^i(\projcomplx) = 0
\]
and
\[
  \H^i(\injcomplx \otimes \projcomplx) = 0
\]
for each $i$.  Using lemma \ref{lem:technical}, this shows
\[
  \{\, \projcomplx \in \KProjR
       \,\mid\, k(\Sigma^i \projcomplx) = 0 \mbox{ for each } i \,\}
  = \EnochsR.
\]
That is, $\EnochsR$ is the kernel of the homological functor $k$.

One consequence of this is that $\EnochsR$ is closed under set indexed
coproducts.  Hence \cite[lem.\ 3.5]{KrauseSmashing} says that for
$\inc$ to have a right-adjoint is the same as for the Verdier quotient
$\KProjR/\EnochsR$ to satisfy that each $\Hom$ set is in fact a set
(as opposed to a class).

Now, the category $\KProjR$ is compactly generated by lemma
\ref{lem:cg}.  By \cite[lem.\ 4.5.13]{NeemanBook} with $\beta =
\aleph_0$, this even implies that there is only a set of isomorphism
classes of compact objects in $\KProjR$.  Hence the version of
Bousfield localization given in \cite[thm.\ 4.1]{PJSpectra} applies to
the functor $k$ on $\KProjR$, and gives that $\KProjR$ modulo the
kernel of $k$ satisfies that each $\Hom$ is a set.  That is,
$\KProjR/\EnochsR$ satisfies that each $\Hom$ is a set, as desired.
\end{proof}

The following elementary result holds by \cite[prop.\
4.10]{PJSpectra}. 

\begin{Lemma}
\label{lem:adjoint_gives_precovering_class}
Let $\inc : \genericsubcategory \longrightarrow \genericcategory$ be
the inclusion of a full subcategory, and suppose that $\adj :
\genericcategory \longrightarrow \genericsubcategory$ is a right
adjoint to $\inc$.  Then $\genericsubcategory$ is a precovering class in
$\genericcategory$. 
\end{Lemma}

Combining proposition \ref{pro:adjoint} and lemma
\ref{lem:adjoint_gives_precovering_class} shows the following.

\begin{Proposition}
\label{pro:Enochs_is_precovering}
The class $\EnochsR$ is precovering in $\KProjR$.
\end{Proposition}

\section{Modules}
\label{sec:modules}

\begin{Lemma}
\label{lem:Gorenstein_projectives_are_precovering}
Let $M$ be an \Rm.  There exists a homomorphism 
\[
  \Gorprojmod \stackrel{g}{\longrightarrow} M 
\]
where $\Gorprojmod$ is a Gorenstein projective \Rm, such that if 
\[
  \Gorprojmodtil \stackrel{\gtil}{\longrightarrow} M 
\]
is any homomorphism with $\Gorprojmodtil$ a Gorenstein projective \Rm\
then there exists a homomorphism
\[
  \Gorprojmodtil \stackrel{\gammatil}{\longrightarrow} \Gorprojmod
\]
so that $\gtil - g\gammatil$ factors through a projective \Rm.
\end{Lemma}

\begin{proof}
Let $\projcomplx \stackrel{\simeq}{\longrightarrow} M$ be a projective
resolution.  Then $\projcomplx$ is in $\KProjR$; let $\Enochscomplx
\stackrel{e}{\longrightarrow} \projcomplx$ be an $\EnochsR$-precover
which exists by proposition \ref{pro:Enochs_is_precovering}.  This gives
\[
  \begin{diagram}[labelstyle=\scriptstyle,height=4ex,width=4ex]
    \cdots & \rTo & \Enochscomplx^{-1} & & \rTo & & \Enochscomplx^0 & & \rTo & & \Enochscomplx^1 & & \rTo & & \Enochscomplx^2 & \rTo & \cdots \\
    & & & & & & & \SE & & \NE & & & & & & & \\
    & & \dTo^{e^{-1}} & & & & \dTo^{e^0} & & \Gorprojmod & & \dTo & & & & \dTo & & \\
    & & & & & & & & \vLine^{g} & & & & & & & & \\
    \cdots & \rTo & \projcomplx^{-1} & & \rTo & & \projcomplx^0 & \rTo & \HonV & & 0 & & \rTo & & 0 & \rTo & \cdots \lefteqn{,} \\
    & & & & & & & \SE & \dTo & & & & & & & & \\
    & & & & & & & & M & & & & & & & & \\
  \end{diagram}
\]
where $\Gorprojmod = \Ker(\Enochscomplx^1 \longrightarrow
\Enochscomplx^2)$ is Gorenstein projective.

Now let $\Gorprojmodtil \stackrel{\gtil}{\longrightarrow} M$ be a
homomorphism with $\Gorprojmodtil$ Gorenstein projective.  Pick
$\Enochscomplxtil$ in $\EnochsR$ so that $\Gorprojmodtil =
\Ker(\Enochscomplxtil^1 \longrightarrow \Enochscomplxtil^2)$.  Clearly
$\gtil$ extends to a chain map $\Enochscomplxtil
\stackrel{\etil}{\longrightarrow} \projcomplx$ so that $\gtil$ and
$\etil$ fit together in
\[
  \begin{diagram}[labelstyle=\scriptstyle,height=4ex,width=4ex]
    \cdots & \rTo & \Enochscomplxtil^{-1} & & \rTo & & \Enochscomplxtil^0 & & \rTo & & \Enochscomplxtil^1 & & \rTo & & \Enochscomplxtil^2 & \rTo & \cdots \\
    & & & & & & & \SE & & \NE & & & & & & & \\
    & & \dTo^{\etil^{-1}} & & & & \dTo^{\etil^0} & & \Gorprojmodtil & & \dTo & & & & \dTo & & \\
    & & & & & & & & \vLine^{\gtil} & & & & & & & & \\
    \cdots & \rTo & \projcomplx^{-1} & & \rTo & & \projcomplx^0 & \rTo & \HonV & & 0 & & \rTo & & 0 & \rTo & \cdots \lefteqn{.} \\
    & & & & & & & \SE & \dTo & & & & & & & & \\
    & & & & & & & & M & & & & & & & & \\
  \end{diagram}
\]

Since $\Enochscomplx \stackrel{e}{\longrightarrow} \projcomplx$ is an
$\EnochsR$-precover, there now exists a chain map $\Enochscomplxtil
\stackrel{\epsilontil}{\longrightarrow} \Enochscomplx$ so that
\begin{equation}
\label{equ:complexes}
  \begin{diagram}[labelstyle=\scriptstyle]
                     &                     & \Enochscomplx \\
                     & \ruTo^{\epsilontil} & \dTo_{e} \\
    \Enochscomplxtil & \rTo_{\etil}        & \projcomplx \\
  \end{diagram}
\end{equation}
is commutative in $\KProjR$.

The chain map $\epsilontil$ induces a homomorphism $\Gorprojmodtil
\stackrel{\gammatil}{\longrightarrow} \Gorprojmod$ so that $\epsilontil$ and
$\gammatil$ fit together in
\[
  \begin{diagram}[labelstyle=\scriptstyle,height=4ex,width=4ex]
    \cdots & \rTo & \Enochscomplxtil^{-1} & & \rTo & & \Enochscomplxtil^0 & & \rTo & & \Enochscomplxtil^1 & & \rTo & & \Enochscomplxtil^2 & \rTo & \cdots \\
    & & & & & & & \SE & & \NE & & & & & & & \\
    & & \dTo^{\epsilontil^{-1}} & & & & \dTo^{\epsilontil^0} & & \Gorprojmodtil & & \dTo_{\epsilontil^1} & & & & \dTo_{\epsilontil^2} & & \\
    & & & & & & & & \vLine^{\gammatil} & & & & & & & & \\
    \cdots & \rTo & \Enochscomplx^{-1} & & \rTo & & \Enochscomplx^0 & \rTo & \HonV & & \Enochscomplx^1 & & \rTo & & \Enochscomplx^2 & \rTo & \cdots \lefteqn{.} \\
    & & & & & & & \SE & \dTo & \NE & & & & & & & \\
    & & & & & & & & \Gorprojmod & & & & & & & & \\
  \end{diagram}
\]
So now there are homomorphisms
\begin{equation}
\label{equ:modules}
  \begin{diagram}[labelstyle=\scriptstyle]
                   &                   & \Gorprojmod \\
                   & \ruTo^{\gammatil} & \dTo_{g} \\
    \Gorprojmodtil & \rTo_{\gtil}      & M \lefteqn{.} \\
  \end{diagram}
\end{equation}

If diagram \eqref{equ:complexes} were commutative as a diagram of
chain maps, then diagram \eqref{equ:modules} would be commutative as a
diagram of modules.  As it is, diagram \eqref{equ:complexes} is only
commutative in $\KProjR$, that is, it is commutative up to chain
homotopy.  It is not hard to see that hence, in diagram
\eqref{equ:modules}, the difference $\gtil - g\gammatil$, while not
necessarily zero, must factor through the module $\Enochscomplxtil^1$.
That is, $\gtil - g\gammatil$ factors through a projective module.
\end{proof}

\begin{Theorem}
\label{thm:Gorenstein_projectives_are_precovering}
Recall setup \ref{set:blanket1}.  In this situation, the Gorenstein
projective modules form a precovering class in the module category of
$\R$. 
\end{Theorem}

\begin{proof}
Let $M$ be a module.  Pick a homomorphism $\Gorprojmod
\stackrel{g}{\longrightarrow} M$ with the property described in 
lemma \ref{lem:Gorenstein_projectives_are_precovering}, and pick a
surjection $\projmod \longrightarrow M$ where $\projmod$ is
projective.  It is easy to see from lemma
\ref{lem:Gorenstein_projectives_are_precovering} that the induced
homomorphism $\Gorprojmod \oplus \projmod \longrightarrow M$ is a
precover with respect to the class of Gorenstein projective modules.
\end{proof}

\section{Non-commutative algebras}
\label{sec:noncomm}

The purpose of this short section is to point out that the above
results apply much more generally than to noetherian commutative rings
with dualizing complexes.  The following setups replace the setups
from the introduction.

\begin{Setup}
\label{set:blanket2}
Let $\R$ be a left-coherent and right-noetherian $k$-algebra over the
field $k$ so that there exists a left-noetherian $k$-algebra
$\matlisR$ and a dualizing complex ${}_{\matlisR}D_{\R}$.  That is,
$D$ is a complex of $\matlisR$-left-$\R$-right-modules, and
\begin{enumerate}

  \item  The cohomology of $D$ is bounded and finitely generated
         both over $\matlisR$ and over $\R^{\opp}$.

  \item  The injective dimensions $\id_{\matlisR} D$ and
         $\id_{\R^{\opp}} D$ are finite.

  \item  The canonical morphisms 
\[
  \R \longrightarrow \RHom_{\matlisR}(D,D) \;\;\; \mbox{and} \;\;\;
  \matlisR \longrightarrow \RHom_{\R^{\opp}}(D,D) 
\]
in the derived categories $\D(\R \otimes_k \R^{\opp})$ and 
$\D(\matlisR \otimes_k \matlisR^{\opp})$ are isomorphisms.

\end{enumerate}
\end{Setup}

\begin{Setup}
\label{set:I2}
Let $D \stackrel{\simeq}{\longrightarrow} \injcomplx$ be an injective
resolution of $D$ over $\matlisR \otimes_k \R^{\opp}$.  Below, I
will replace $\injcomplx$ by a bounded truncation.  This may ruin the
property that $\injcomplx$ is an injective resolution over
$\matlisR \otimes_k \R^{\opp}$, but because $\id_{\matlisR}
D$ and $\id_{\R^{\opp}} D$ are finite, I can still suppose that
$\injcomplx$ consists of modules which are injective both over
$\matlisR$ and over $\R^{\opp}$.
\end{Setup}

The above definition of dualizing complexes over non-commutative
algebras is due to \cite[def.\ 1.1]{YekutieliZhang}.

With setups \ref{set:blanket1} and \ref{set:I1} replaced by setups
\ref{set:blanket2} and \ref{set:I2}, let me inspect the rest of the
paper.  As the ground ring $\R$ is now non-commutative, I must replace
``module'' by ``left-module'' throughout.  Remark
\ref{rmk:flats_have_finite_pd} also needs to be replaced by the
following.

\begin{Remark}
\label{rmk:flats_have_finite_pd2}
The results of \cite{PJfdpd} apply under setup \ref{set:blanket2}, and
so each flat \Rlm\ has finite projective dimension.
\end{Remark}

After this, the proof of lemma \ref{lem:technical} goes through if one
keeps track of left and right structures throughout.  The proofs of
lemma \ref{lem:cg} and proposition \ref{pro:adjoint} also still work,
and along with lemma \ref{lem:adjoint_gives_precovering_class} this
still implies proposition \ref{pro:Enochs_is_precovering}.  And
finally, the proofs of lemma
\ref{lem:Gorenstein_projectives_are_precovering} and theorem
\ref{thm:Gorenstein_projectives_are_precovering} still go through.

So theorem \ref{thm:Gorenstein_projectives_are_precovering} remains
valid.  Let me formulate this in full.

\begin{Theorem}
\label{thm:Gorenstein_projectives_are_precovering2}
Recall setup \ref{set:blanket2}.  In this situation, the Gorenstein
projective modules form a precovering class in the category of \Rlm s.
\end{Theorem}

\begin{Corollary}
\label{cor:PI}
Let $\generalR$ be a noetherian $k$-algebra and suppose that one of the
following holds.
\begin{enumerate}

  \item  $\generalR$ is a complete semi-local PI algebra.

  \item  $\generalR$ has a filtration $F$ so that the associated
         graded algebra $\operatorname{gr}^F\!\generalR$ is connected
         and noetherian, and either PI, graded FBN, or with enough
         normal elements. 

\end{enumerate}
Then the Gorenstein projective modules form a precovering class in
the category of \Rlm s.
\end{Corollary}

\begin{proof}
The algebra $\generalR$ can be used as $\R$ in setup
\ref{set:blanket2} and theorem
\ref{thm:Gorenstein_projectives_are_precovering2} because $\matlisR$
and $D$ exist.  In case (i) this is by \cite[cor.\
0.2]{WuZhangDualizing}, and in case (ii) by \cite[cor.\
6.9]{YekutieliZhang}.
\end{proof}

\bigskip

\noindent
{\bf Acknowledgement.}  
I thank Juan Antonio L\'{o}pez-Ramos for communicating
\cite{EnochsJendaLopezRamos} and \cite{EnochsLopezRamos}.

The diagrams were typeset with Paul Taylor's
{\tt diagrams.tex}.

\end{document}